\documentclass[aop,seceqn,MSNbibl,citesort,dvips]{arximspdf}

\doi{10.1214/10-AOP576}
\volume{39}
\issue{3}
\pubyear{2011}
\firstpage{1137}
\lastpage{1160}

\makeatletter

\newtheorem{proposition}{Proposition}
\newtheorem{theorem}{Theorem}
\newproclaim{example}{Example}
\newproclaim{definition}{Definition}
\newproclaim{remark}{Remark}

\newcommand{\R}{\mathbb{R}}
\newcommand{\bx}{\mathbf{X}}
\newcommand{\F}{\mathcal{F}}
\newcommand{\G}{\mathcal{G}}

\makeatother

\begin{document}
\begin{frontmatter}

\title{Nonparametric sequential prediction for stationary~processes}
\runtitle{Nonparametric sequential prediction}

\begin{aug}
\author[A]{\fnms{Guszt\'av} \snm{Morvai}\thanksref{t1}\ead
[label=e1]{morvai@math.bme.hu}\ead
[label=u1,url]{http://www.math.bme.hu/\textasciitilde morvai/}} and
\author[B]{\fnms{Benjamin} \snm{Weiss}\corref{}\ead[label=e2]{weiss@math.huji.ac.il}}
\runauthor{G. Morvai and B. Weiss}
\affiliation{MTA-BME Stochastics Research Group and Hebrew University
of Jerusalem}
\address[A]{MTA-BME Stochastics Research Group\\
Institute of Mathematics\\
Egry J\'ozsef utca 1\\
Building H \\
Budapest, 1111\\
Hungary\\
\printead{e1}\\
\printead{u1}}
\address[B]{Hebrew University of Jerusalem\\
Institute of Mathematics\\
Jerusalem 91904\\
Israel\\
\printead{e2}}
\end{aug}

\thankstext{t1}{Supported by the Bolyai J\'anos Research Scholarship
and OTKA Grant K75143.}

\received{\smonth{11} \syear{2007}}
\revised{\smonth{5} \syear{2010}}

%
\begin{abstract}
We study the problem of finding an universal estimation scheme
$h_n\dvtx \R^n \to\R$, $n=1,2,\ldots$ which will satisfy
\begin{eqnarray*}
&&\lim_{t\rightarrow\infty} {1\over t} \sum_{i=1}^{t}
|h_i(X_0,X_1,\ldots,X_{i-1})\\
&&\qquad\hspace*{24.75pt}{} -
E(X_{i}|X_0,X_1,\ldots,X_{i-1})|^p=0  \qquad\mbox{a.s.}
\end{eqnarray*}
for all real valued stationary and ergodic processes that are in $L^p$.
We will construct a single such scheme for all $1 < p \le\infty$, and
show that for $p=1$ mere integrability does not suffice but $L\log^+L$
does.
\end{abstract}

%

\end{frontmatter}

\section{Introduction}
The problem of sequentially predicting the next value $X_n$ of
a stationary process after observing the initial values $X_i$ for
$0 \le i < n$ is one of the central problems in probability and statistics.
Usually, one bases the prediction on the conditional expectation
$E(X_n|X_0^{n-1})$
where we write for brevity $X_0^{n-1}=\{X_0,X_1,\ldots,X_{n-1}\}$. However, when
one does not know the distribution of the process one is faced with the
problem of estimating the conditional expectation from a single sample
of length
$n$. It was shown long ago by Bailey \cite{Bailey1976} (cf. also Ryabko
\cite{Ryabko88} and
Gy\"orfi Morvai and Yakowitz \cite{GYMY98}) that even
for binary
processes no universal scheme $h_n(X_0^{n-1})$ exists which will almost surely
satisfy $ \lim_{n\to\infty} (h_n(X_0^{n-1}) - E(X_n|X_0^{n-1})) =0$.
This is in
contrast
to the backward estimation problem where one is trying to estimate
$E(X_0|X^{-1}_{-\infty})$ based on the successive observations of
$X^{-1}_{-\infty}$.
Here, it was Ornstein \cite{Ornstein1978} who constructed the first such
universal estimator for finite valued processes. This was generalized
to bounded
processes by Algoet \cite{Algoet1992}, Morvai \cite{MorvaiPhD} and Morvai
Yakowitz and Gy\"orfi~\cite{MYGY1996}. For unbounded processes,
several universal
estimators were constructed (see Algoet
\cite{Algoet1999} and Gy\"orfi et al. \cite{GKKW2002}).

Returning to our original problem of sequential prediction
it was already observed by Bailey that backward schemes could be
used for the sequential prediction problem successfully in the
sense that that the error tends to zero in the Ces\'aro mean. To
establish this, he applied a generalized ergodic theorem which
requires some technical hypotheses which were satisfied in his
case.

Over the years some authors have extended this work, namely
of adapting backward schemes to sequential prediction, but only for bounded
processes (see Algoet \cite{Algoet1992,Algoet1999}, Morvai \cite
{MorvaiPhD}, Morvai
Yakowitz and Gy\"orfi \cite{MYGY1996}
and Gy\"orfi et al. \cite{GKKW2002}).

Another approach to the sequential prediction used a weighted average of
expert schemes, and with these results were extended to the general unbounded
case by Nobel \cite{Nobel2003}
and Ottucsak \cite{GyO2007} (see also the survey of Feder and Merhav
\cite{FM98}).
However, none of these results were optimal in the sense that moment conditions
higher than necessary were assumed.
It is our
purpose to obtain these optimal conditions and to show why they are necessary.
We consider the following problem for $1 \le p \le\infty$. Does there exist a scheme
$h_n(X_0^{n-1})$ which will satisfy
\[
\lim_{t\rightarrow\infty} {1\over t} \sum_{i=1}^{t} |h_i(X_0^{i-1})-
E(X_{i}|X_{0}^{i-1})|^p=0  \qquad\mbox{a.s.}
\]
for all real valued stationary and ergodic processes that are in $L^p$.
The only
case that has been solved completely is when $p$ is infinity. Even the recent
schemes Nobel \cite{Nobel2003} and Gy\"orfi
and Ottucsak \cite{GyO2007} put a higher moment
condition
on the process than is manifestly required. Our main result is that the basic
scheme first introduced by the first author in his thesis can be
adapted to give
a scheme which will answer our problem positively for all $1 < p$. For
$p=1$, we shall
show that stronger hypothesis is necessary, as is usually the case, and
we will
establish the convergence under the hypothesis that $X_0 \in L\log^+L$.

In the third section, we will show how this hypothesis cannot be
weakened to
$X_0 \in L^1$. Our construction will be based on one of the simplest ergodic
transformation, the adding machine, and illustrates the richness of behavior
that is possible for processes that are almost periodic (in the
sense of Besicovich).

As soon as one knows that the errors converge to zero in Ces\'aro mean,
it follows that there is a set of density one of time moments along
which the errors converge to zero. However, in general one does not know
what this sequence is. In the framework of estimation, schemes adapted
to a sequence of stopping times (see
\cite{Mo00,MW03,MW05ThSP,MW05StatLett,MW05PTRF,MW04TEST,MW05Poinc2})
one may ask can one find a sequence of stopping times with density one
along which the errors of a universal sequential prediction scheme will
tend to zero. We have been unable to do this in general and regard it
as an important open problem. Finally, we refer the interested reader
to some other papers which are relevant to this line of research
\cite{Algoet1994,GyLM1999,MYA1997,MW2005Poinc1,Weiss00}.

Some technical probabilistic results have been relegated to the
\hyperref[app]{Appendix}, they are
of a classical nature and may be known,
but we were unable to find references.

\section{The main result}

Let $\bx=\{X_n\}$ denote a real-valued doubly infinite
stationary ergodic time series. Let
\[
X^{j}_{i}=(X_{i},X_{i+1},\ldots,X_{j})
\]
be notation for a data
segment, where $i$ may be minus infinity. Let
\[
\bx
^-=X_{-\infty}^{-1}.
\]

Let $G_k$ denote the quantizer
\[
G_k(x)=
\cases{
0, &\quad if $-2^{-k}<x<2^{-k}$,\cr
-i2^{-k},  &\quad if $-(i+1)2^{-k}<
x\le-i2^{-k}$ for some $i=1,2,\ldots,$\cr
i2^{-k}, &\quad if $i2^{-k}\le x<(i+1)2^{-k}$.}
\]
Define the sequences $\lambda_{k-1}$ and $\tau_k$ recursively
($k=1,2,\ldots$).
Put $\lambda_0=1$ and let $\tau_k$ be the time between the
occurrence of the pattern
\[
B(k)=(G_k(X_{-\lambda_{k-1}}),\ldots,G_k(X_{-1}))=
G_k(X_{-\lambda_{k-1}}^{-1})
\]
at time $-1$ and the last occurrence of the same pattern prior to time
$-1$. More precisely, let
\[
{\tau}_k=
\min\{t>0\dvtx G_k(X_{-\lambda_{k-1}-t}^{-1-t})=G_k(X_{-\lambda
_{k-1}}^{-1})\}.
\]
Put
\[
\lambda_k=\tau_k+\lambda_{k-1}.
\]
Define
%
\begin{equation}
\label{eqalgregk}
R_k={1\over k}\sum_{1\le j\le k} X_{-\tau_j}.
\end{equation}
To obtain a fixed sample size $t>0$ version,
let $\kappa_t$ be the maximum of integers $k$ for which $\lambda_k
\le t$.
For $t>0$, put
%
\begin{equation}
\label{eqalgregt}
\hat R_{-t}={1\over\kappa_t}\sum_{1\le j\le\kappa_t} X_{-\tau_j}.
\end{equation}
Motivated by Bailey \cite{Bailey1976}, for $t>0$ consider the estimator
\[
\hat R_{t}(\omega)=\hat R_{-t}(T^t\omega),
\]
which is defined in terms of
$(X_0,\ldots,X_{t-1})$ in the same way as $\hat R_{-t}(\omega)$ was defined
in terms of $(X_{-t},\ldots,X_{-1})$. ($T$ denotes the left shift
operator.)
The estimator $\hat R_t$ may be viewed as an online predictor of $X_t$.
This predictor has special significance not only because of potential
applications, but additionally because
Bailey \cite{Bailey1976} proved that it is impossible to construct
estimators $\hat
R_t$
such that always $\hat R_t-E(X_t|X^{t-1}_0)\to0$ almost surely.
\begin{theorem} \label{positivetheorem}
Let $\{X_n\}$ be stationary and ergodic. Assume that
\[
E(|X_0|\log^+(|X_0|))<\infty.
\]
Then
%
\begin{equation}
\label{eqregplust1}
\lim_{t\rightarrow\infty} {1\over t} \sum_{i=1}^{t} |\hat R_{i}-
E(X_{i}|X_{0}^{i-1})|=0  \qquad\mbox{a.s.}
\end{equation}
and
%
\begin{equation}
\label{eqregplust1.2}
\lim_{t\rightarrow\infty} {1\over t} \sum_{i=1}^{t} |\hat R_{i}-
X_i|=E\bigl(|E(X_{0}|X_{-\infty}^{-1})-X_0|
\bigr)  \qquad\mbox{a.s.}
\end{equation}
Furthermore, if for some $1<p<\infty$, $E(|X_0|^p)<\infty$, then
%
\begin{equation}
\label{eqregcesaro}
\lim_{t\rightarrow\infty} {1\over t} \sum_{i=1}^{t} |\hat R_{i}-
E(X_{i}|X_{0}^{i-1})|^p=0  \qquad\mbox{a.s.}
\end{equation}
and
%
\begin{equation}
\label{eqregcesaro1.2}
\lim_{t\rightarrow\infty} {1\over t} \sum_{i=1}^{t} |\hat R_{i}-
X_{i}|^p=E\bigl(|E(X_{0}|X_{-\infty}^{-1})-X_0
|^p\bigr)  \qquad\mbox{a.s.}
\end{equation}
\end{theorem}
\begin{pf}
The proof will follow the same pattern in all four cases. We will verify
that the backward estimator scheme converges almost surely and we will
see that the sequence of errors is dominated by an integrable function.
This allows us to conclude from
the generalized ergodic theorem of Maker (rediscovered by Breiman,
cf. Theorem 1 in Maker \cite{Maker1940} or Theorem 12 in Algoet
\cite{Algoet1994}) that the forward scheme converges in Cesaro mean.
For the first case, we will carry this out in full detail, for the
others we will just check the requisite properties for the
backward scheme. First, consider
\begin{eqnarray*}
R_k &=& {1\over k} \sum_{1\le j\le k}
[X_{-\tau_j}-G_j(X_{-\tau_j})] \\
&&{}+ {1\over k} \sum_{1\le j\le k}
[ G_j(X_{-\tau_j})
- E(G_j(X_{-\tau_j})|G_{j-1}(X_{-\lambda_{j-1}}^{-1}))]\\
&&{}+ {1\over k}\sum_{1\le j\le k}
[ E(G_j(X_{-\tau_j})|G_{j-1}(X_{-\lambda_{j-1}}^{-1}))-
E(X_{-\tau_j}|G_{j-1}(X_{-\lambda_{j-1}}^{-1}))]\\
&&{}+ {1\over k}\sum_{1\le j\le k}
[ E(X_{-\tau_j}|G_{j-1}(X_{-\lambda_{j-1}}^{-1}))-
E(X_{0}|G_{j-1}(X_{-\lambda_{j-1}}^{-1}))]\\
&&{}+ {1\over k}\sum_{1\le j\le k}
E(X_{0}|G_{j-1}(X_{-\lambda_{j-1}}^{-1}))\\
&=&A_k+B_k+C_k+D_k+E_k.
\end{eqnarray*}
Obviously,
\[
|A_k|+|C_k|\le{2\over k}\sum_{1\le j\le k} 2^{-j}\le{2\over k}\to0.
\]
Now we will deal with $D_k$.
By Lemma 1, in Morvai, Yakowitz and Gy\"orfi \cite{MYGY1996},
\[
P\bigl(X_{-\tau_j}\in C|G_{j-1}(X_{-\lambda_{j-1}}^{-1})\bigr)=
P\bigl(X_0\in C|G_{j-1}(X_{-\lambda_{j-1}}^{-1})\bigr).
\]
Using this, we get that $D_k=0$.

Assume that $E(|X_0|\log^+(|X_0|))<\infty$.
Toward mastering $B_k$, one observes that
$\{X_{-\tau_j}\}$
are identically distributed
by Lemma 1 in Morvai, Yakowitz and Gy\"orfi \cite{MYGY1996}
and
$B_k$ is an average of
martingale differences.
By Proposition \ref{PropositionElton} in the
\hyperref[app]{Appendix}, $|B_k|\to0$ almost surely and $E(\sup_{1\le k}
|B_k|)<\infty$.

Now we deal with the last term $E_k$.
By assumption,
\[
\sigma(G_{j}(X_{-\lambda_{j}}^{-1}))\uparrow\sigma(\bx^-).
\]
Consequently by the a.s. martingale convergence theorem, we have that
\[
E(X_0|G_{j}(X_{-\lambda_{j}}^{-1}))\to E(X_0|\bx^-)
\qquad\mbox{a.s.,}
\]
and thus
\[
E_k\to E(X_0|\bx^-) \qquad\mbox{a.s.}
\]
Furthermore,
by Doob's inequality, cf.
Theorem 1 on page 464, Section 3, Chapter~VII in Shiryayev \cite{Shiryayev1984},
$E({\sup_{1\le k}} |E_k|)\le E({\sup_{1\le j}} |E
(X_0|G_{j}(X_{-\lambda_{j}}^{-1}))|)<\infty$.

We have so far proved that
\[
R_k \rightarrow E(X_0|\bx^-) \qquad\mbox{almost surely}
\]
and
\[
E \Bigl( {\sup_{1\le k}} |R_k|\Bigr) <\infty.
\]
This in turn implies that
\[
\lim_{t\to\infty}{\hat R}_{-t}= E(X_0|\bx^-) \qquad\mbox{almost surely}
\]
and
\[
E \Bigl( {\sup_{1\le t}} |{\hat R}_{-t}|\Bigr) <\infty.
\]
Now since $E(X_0|X^{-1}_{-t})\to E(X_0|\bx^-)$ almost surely,
\[
{\lim_{t\to\infty}}|{\hat R}_{-t}-E(X_0|X^{-1}_{-t})|=0
\qquad\mbox{almost surely}
\]
and
by Doob's inequality,
%
\begin{eqnarray*}
E \Bigl( {\sup_{1\le t}} |{\hat R}_{-t}-E(X_0|X^{-1}_{-t})|\Bigr)
&\le& E \Bigl( {\sup_{1\le t}} |{\hat R}_{-t}|\Bigr)\\
&&{}+ E \Bigl( {\sup_{1\le t}} |E(X_0|X^{-1}_{-t})|\Bigr)\\
&<& \infty.
\end{eqnarray*}
Now, apply the generalized ergodic theorem
to conclude that
\begin{eqnarray*}
\lim_{t\rightarrow\infty} {1\over t} \sum_{i=1}^{t} \bigl( |
{\hat R}_{-i}-
E(X_{0}|X_{-i}^{-1})|(T^i\omega)\bigr) &=&
\lim_{t\rightarrow\infty} {1\over t} \sum_{i=1}^{t} |\hat R_{i}-
E(X_{i}|X_{0}^{i-1})|\\
&=&0  \qquad\mbox{a.s.}
\end{eqnarray*}
and the proof of (\ref{eqregplust1}) is complete. Similarly,
\[
{\lim_{t\to\infty}} |{\hat R}_{-t}-X_0|=|E(X_0|X^{-1}_{-\infty
})-X_0|
\qquad\mbox{almost surely}
\]
and
\[
E \Bigl( {\sup_{1\le t}} |{\hat R}_{-t}-X_0|\Bigr)
\le E \Bigl( {\sup_{1\le t}} |{\hat R}_{-t}|\Bigr)
+ E ( |X_0|) < \infty
\]
and the generalized ergodic theorem gives
\begin{eqnarray*}
\lim_{t\rightarrow\infty} {1\over t} \sum_{i=1}^{t} \bigl( |
{\hat R}_{-i}-
X_{0}|(T^i\omega)\bigr) &=&
\lim_{t\rightarrow\infty} {1\over t} \sum_{i=1}^{t} |\hat R_{i}-
X_{i}|\\
&=&
E\bigl(|E(X_0|X^{-1}_{-\infty})-X_0|\bigr)
\qquad\mbox{a.s.}
\end{eqnarray*}
and the proof of (\ref{eqregplust1.2}) is complete.

Now, we assume that for some $1<p<\infty$, $E(|X_0|^p)<\infty$, and
we prove
(\ref{eqregcesaro}).
Observe that
\[
|R_k|^p \le3^p \biggl[ \biggl({2\over k}\biggr)^p
+|B_k|^p+|E_k|^p\biggr]
\]
and since by Proposition \ref{PropositionL}
in the \hyperref[app]{Appendix}
$|B_k|\to0$ almost surely and $E({\sup_{1\le k}}
|B_k|^p)<\infty$
and
by Doob's inequality,\vspace*{1pt}
$E({\sup_{1\le k}} |E_k|^p)<\infty$ and $E_k\to E(X_0|\bx^-)$ almost surely
(for the same
reason as before).

We have so far proved that
\[
R_k \rightarrow E(X_0|\bx^-) \qquad\mbox{almost surely}
\]
and
\[
E \Bigl( {\sup_{1\le k} }|R_k|^p\Bigr) <\infty.
\]
This in turn implies that
\[
\lim_{t\to\infty}{\hat R}_{-t}= E(X_0|\bx^-) \qquad\mbox{almost surely}
\]
and
\[
E \Bigl( {\sup_{1\le t} }|{\hat R}_{-t}|^p\Bigr) <\infty.
\]
Now since $E(X_0|X^{-1}_{-t})\to E(X_0|\bx^-)$ almost surely,
\[
{\lim_{t\to\infty}}|{\hat R}_{-t}-E(X_0|X^{-1}_{-t})|^p=0
\qquad\mbox{almost surely}
\]
and
by Doob's inequality,
%
\begin{eqnarray*}
E \Bigl( {\sup_{1\le t} }|{\hat R}_{-t}-E(X_0|X^{-1}_{-t})|^p\Bigr)
&\le& 2^p E \Bigl( {\sup_{1\le t}} |{\hat R}_{-t}|^p\Bigr)\\
&&{}+ 2^p E \Bigl( {\sup_{1\le t}} |E(X_0|X^{-1}_{-t})|^p\Bigr)\\
&<& \infty.
\end{eqnarray*}
By
Maker's (or Breiman's) generalized ergodic theorem (cf. Theorem 1
in Maker \cite{Maker1940} or Theorem 12 in Algoet
\cite{Algoet1994})
one gets
(\ref{eqregcesaro}).
Similarly,
\[
{\lim_{t\to\infty}} |{\hat R}_{-t}-X_0|^p=|E(X_0|X^{-1}_{-\infty
})-X_0|^p
\qquad\mbox{almost surely}
\]
and
\[
E \Bigl( {\sup_{1\le t} }|{\hat R}_{-t}-X_0|^p\Bigr)
\le2^p E \Bigl( {\sup_{1\le t}} |{\hat R}_{-t}|^p\Bigr)
+ 2^p E ( |X_0|^p) < \infty.
\]
Now, apply Maker's (or Breiman's) generalized ergodic theorem
to prove
(\ref{eqregcesaro1.2}).
The proof of Theorem \ref{positivetheorem} is complete.
\end{pf}
\begin{remark}
We are indebted to the referee for the following remark. Using the
notion of
Bochner integrability of strongly measurable functions with values in $c_0$
and the extension of Birkhoff's ergodic theorem to Banach space valued functions
(see Krengel \cite{Krengel1985}, page 167),
one can give an easy proof of Maker's theorem. The key condition now
becomes the
fact that the norm of the sequence $\{f-f_k\}$ in $c_0$ is integrable,
and then the
convergence in the norm of $c_0$ allows one to deduce the convergence
of the
diagonal sequence which is what appears in Maker's theorem.
\end{remark}

\section{Integrability alone is not enough}

In Theorem \ref{positivetheorem} for the Ces\'aro convergence in the $L^1$
norm, we assumed that $X_0$ was not merely in $L^1$ but in
$L^1\log^+L$. In this section, we shall show that some additional
condition is really necessary. We will first give an example to show
that the maximal function of the conditional expectations
${\sup_{1 \le n}} |E(X_0|X^{-1}_{-n})|$ may be nonintegrable for an integrable
process. We shall do so in an indirect fashion by showing that the
the estimate $E(X_n|X^{n-1}_0)$ for $E(X_n|X^{n-1}_{-\infty})$ does not
converge in Ces\'aro mean to zero. This means that even though we are
may be
in the distant future the information of the prehistory can make a serious
difference. This example serves as a model for the main result of the section
where we show that for any estimation scheme for $E(X_n|X^{n-1}_0)$ which
converges almost surely in Ces\'aro mean for all bounded processes there
will be some ergodic integrable process where it fails to converge. Indeed
the processes that we need to consider are countably valued and in fact
are zero entropy and finitarily Markovian (see below for a definition),
a generalization of finite order Markov chains.

First, let us fix the notation.
Let $\{X_n\}_{n=-\infty}^{\infty}$ be a stationary and ergodic time
series taking
values from finite or countable alphabet
$\mathcal{X}$. (Note that all stationary time series $\{X_n\}
_{n=0}^{\infty}$
can be thought to be a
two sided time series, that is, $\{X_n\}_{n=-\infty}^{\infty}$.)
\begin{definition}
The stationary time series $\{X_n\}$ is said to be finitarily Markovian if
almost surely the sequence of the conditional distributions $\mathcal
L$ $(X_1|X^0_{-k})$
is constant for large $k$ (it is random how large $k$ should be).
\end{definition}

This class includes of course all finite order Markov chains but also many
other processes such as the finitarily determined processes of Kalikow,
Katznelson and Weiss \cite{KKW92},
which serve to represent all isomorphism classes of zero entropy processes.

For some concrete examples that are not Markovian, consider the following
example.
\begin{example}
Let $\{M_n\}$ be any stationary and ergodic first order Markov chain with
finite or countably infinite
state space $S$.
Let $s\in S$ be an arbitrary state with $P(M_1=s)>0$. Now let $X_n=I_{\{
M_n=s\}}$.
By Shields \cite{Sh96}, Chapter I.2.c.1, the binary time series $\{X_n\}$
is stationary and ergodic.
It is also finitarily Markovian. Indeed, the conditional probability
$P(X_1=1|X^0_{-\infty})$
does not depend on values
beyond the first (going backward) occurrence of one in $X^0_{-\infty}$
which identifies the first (going backward) occurrence of state $s$ in the
Markov chain $\{M_n\}$.
The resulting time series $\{X_n\}$ is not a Markov chain of any order
in general.
\end{example}

We note that Morvai and Weiss \cite{MW05Bern} proved that there is no
classification
rule
for discriminating the class of finitarily Markovian processes from
other ergodic
processes. For more about estimation for finitarily Markovian
processes, see
Morvai and Weiss \cite{MW05StatLett,MW05PTRF,MW05Poinc2}.
\begin{theorem} \label{counterex1}
Let $\mathcal{X}=\{0, 10^{-k}, {2^{k}\over3^{m}},
k=1,2,\ldots, m=1,2,\ldots\}$.
There exists a stationary and ergodic finitarily Markovian time series
$\{X_n\}$
taking values from $\mathcal{X}$ such that $E |X_0|<\infty$
and
\[
\limsup_{N\to\infty} {1\over N} \sum_{n=1}^{N} |E
(X_n|X^{n-1}_0)-
E(X_n|X^{n-1}_{-\infty})|=\infty
\]
almost surely. Therefore,
\[
E\Bigl({\sup_{1\le n} }|E(X_0|X^{-1}_{-n})
|\Bigr)=\infty.
\]
\end{theorem}
\begin{pf}
Let $\Omega$ be the one sided sequence space over $\{0,1\}$. Let
$\omega=(\omega_1,\omega_2,\ldots)\in\Omega$. Define the
transformation $T\dvtx
\Omega\rightarrow\Omega$ as follows:
\[
(T\omega)_i=
\cases{
0,  &\quad if $\omega_j=1$ for all $j\le i$,\cr
1,  &\quad if $\omega_i=0$ and for all $j<i\dvtx\omega_j=1$,\cr
\omega_i,  &\quad otherwise.}
\]
Consider the product measure $P=\Pi_{i=1}^{\infty} \{1/2,1/2\}$ on
$\Omega$
which is preserved by $T$. It is well known (cf. Aaronson
\cite{Aaronson1997}, page 25) that
$(\Omega, P,T)$ is an ergodic
process, called the adding machine or dyadic odometer.
The process will be defined by a function
$f\dvtx \Omega\rightarrow\R$ as $X_n(\omega)=f(T^n\omega)$.
Let $l_3<\cdots< l_{k-1}< l_k\to\infty$.
Define $a_k=a$ and $b_k=b$ when $k=2^a+b$ where $1\le b\le2^a$.
Define
\[
C_k=\{\omega\dvtx \omega_i=1  \mbox{ for }1\le i<l_k, \omega_{l_k}=0\},
\]
clearly $P(C_k)=2^{-l_k}$. Let
\[
D_k=\{\omega\dvtx \omega_i=1  \mbox{ for }1\le i<l_k-a_k, \omega_{l_k-a_k}=0,
\omega_i=1  \mbox{ for }l_k-a_k< i<l_k\}
\]
and
\[
E_k=\bigcup_{i=0}^{2^{l_k-a_k-1}-1}T^{-i} D_k.
\]
Notice that
\[
E_k=\{\omega\dvtx \omega_{l_k-a_k}=0, \omega_{j}=1\mbox{,
for all }l_k-a_k< j< l_k\}.
\]
It is clear that if the $l_k$'s are chosen large enough so that for all $k'>k$
$l_k<l_{k'}-2 a_{k'}$:
\begin{itemize}
\item the family $C_k,D_l$ $k,l\ge3$ consists of disjoint sets,
\item the intervals $[l_k-a_k,l_k-1]$ are also disjoint and therefore
the sets
$E_k$ are independent.
\end{itemize}
The signaling function $u$ is defined by
\[
u(\omega)=\sum_{k=3}^{\infty} 10^{-k} I_{D_k}(\omega)
\]
and the main contributor to $f$ will be
\[
v(\omega)=\sum_{k=3}^{\infty} {2^{l_k} \over3^{a_k} }
I_{C_k}(\omega).
\]
Clearly,
\begin{eqnarray*}
E(v(\omega))&=&\sum_{k=3}^{\infty} {2^{l_k} \over3^{a_k} } P(C_k)=
\sum_{k=3}^{\infty} {1 \over3^{a_k} } =
\sum_{a=1}^{\infty} \sum_{b=1}^{2^a} {1 \over3^{a} }\\
&=&
\sum_{a=1}^{\infty} \biggl({2 \over3}\biggr)^{a}<\infty.
\end{eqnarray*}
Define a process by
$
f(\omega)=u(\omega)+v(\omega)
$
and
\[
X_n(\omega)=f(T^n\omega).
\]
Notice that $X_n\in\{0, 10^{-k}, {2^{l_k}\over3^{a_k}} ,
k=3,4,\ldots\}$.
Observe that $P(E_k)=2^{-a_k}$ and
\[
\sum_{k=3}^{\infty} P(E_k)=\sum_{a=1}^{\infty} \sum_{b=1}^{2^a} 2^{-a}=
\sum_{a=1}^{\infty} 1= \infty.
\]
By the Borel--Cantelli lemma, a point $\omega$ belongs to $E_k$
infinitely often.
When \mbox{$\omega\in E_k$},
\[
T^{i_0}\omega\in D_k \qquad\mbox{for some }0\le i_0\le2^{l_k-a_k-1}-1.
\]
For $\omega\in E_k$, we know that $X_{i_0}(\omega)=10^{-k}$.
At time $i_0+2^{l_k-a_k-1}-1$,
\[
(T^{i_0+2^{l_k-a_k-1}-1}(\omega))_j =
\cases{
0,  &\quad if $j=1$,\cr
1,  &\quad if $1<j\le l_k-1$,\cr
\omega_j, &\quad otherwise.}
\]
Let's compute for a fixed $i_0$ such that $T^{i_0}\omega\in D_k$ (i.e.,
$X_{i_0}=10^{-k}$)
\[
E( X_{i_0+2^{l_k-a_k}}| X^{i_0+2^{l_k-a_k}-1}_0).
\]
Take $N=2^{l_k-a_k}$ and consider
\[
{1\over N} \sum_{n=1}^{N} | E( X_{n}| X^{n-1}_0)-
E( X_{n}|
X^{n-1}_{-\infty})|.
\]
For $\omega\in T^{-i_0}D_k$ (i.e., $X_{i_0}=10^{-k}$ ), we know that
\[
(T^{i_0+2^{l_k-a_k-1}}\omega)_j =
\cases{
1,  &\quad if $1\le j \le l_k-1$,\cr
\omega_j, &\quad otherwise.}\vadjust{\goodbreak}
\]
Therefore if $X_{i_0+2^{l_k-a_k-1}}>0$, then we must have
\[
T^{i_0+2^{l_k-a_k-1}} \omega\in C_k\cup\bigcup_{j>k} (C_j\cup D_j)
\]
(because if $k'<k$ then $l_{k'}<l_k$ and $C_{k'}$, $D_{k'}$ are defined
by zero
values of $\omega_i$ with $i<l_k$)
and
\begin{eqnarray*}
&&E( X_{i_0+2^{l_k-a_k-1}}| X^{i_0+2^{l_k-a_k-1}-1}_0)\\
&&\qquad=
{ {2^{l_k} /3^{a_k}} 2^{-l_k} +\sum_{j>k} {2^{l_j} / 3^{a_j}}
2^{-l_j}+
\sum_{j>k} 10^{-j} 2^{-l_j} +0 \over
P(D_k) } \\
&&\qquad\ge {({2/3})^{a_k+1} \over2\cdot2^{-l_k} }\\
&&\qquad={1\over2} 2^{l_k} \biggl({2\over3}\biggr)^{a_k+1}.
\end{eqnarray*}
Similarly,
\begin{eqnarray*}
&&E( X_{i_0+2^{l_k-a_k}}| X^{i_0+2^{l_k-a_k}-1}_0)\\
&&\qquad=
{ {2^{l_k} /3^{a_k}} 2^{-l_k} +\sum_{j>k} {2^{l_j} /3^{a_j}}
2^{-l_j}+
\sum_{j>k} 10^{-j} 2^{-l_j} +0 \over
P(D_k) } \\
&&\qquad\le {10^{-k-1}+\sum_{i=0}^{\infty} ({2/3}
)^{a_k+i} \over2\cdot
2^{-l_k} }\\
&&\qquad=
{1\over2} 2^{l_k} \biggl( 10^{-k-1}+\biggl({2\over3}\biggr)^{a_k}
3\biggr)\\
&&\qquad\le 4 \cdot2^{l_k} \biggl({2\over3}\biggr)^{a_k}.
\end{eqnarray*}
On the other hand,
$X^{i_0+2^{l_k-a_k-1}-1}_{-\infty}$ determines exactly the value of
$X_{i_0+2^{l_k-a_k-1}}$.
There are four cases. If $X_{i_0+2^{l_k-a_k-1}}$ is equal with
$0$, ${2^{l_k}\over3^{a_k}}$, or for some $k<k'\dvtx 10^{-k'}$ or
${2^{l_{k'}}\over3^{a_{k'}}}$. That is,
\begin{eqnarray*}
&&E( X_{i_0+2^{l_k-a_k-1}}| X^{i_0+2^{l_k-a_k-1}-1}_{-\infty
})\\
&&\qquad=
\cases{
\displaystyle {2^{l_k}\over3^{a_k}},  &\quad if $T^{i_0+2^{l_k-a_k-1}}\omega\in
C_k$,\vspace*{2pt}\cr
\displaystyle 10^{-k'}, &\quad if $T^{i_0+2^{l_k-a_k-1}}\omega\in D_{k'}$ for some
$k<k'$,\vspace*{2pt}\cr
\displaystyle {2^{l_{k'}}\over3^{a_{k'}}},  &\quad if $T^{i_0+2^{l_k-a_k-1}}\omega
\in C_{k'}$ for some $k<k'$,\vspace*{2pt}\cr
0,  &\quad if otherwise.}
\end{eqnarray*}
Now
\begin{eqnarray*}
&&
| E( X_{i_0+2^{l_k-a_k-1}}| X^{i_0+2^{l_k-a_k-1}-1}_0
)- E(
X_{i_0+2^{l_k-a_k-1}}| X^{i_0+2^{l_k-a_k-1}-1}_{-\infty})|
\\
&&\qquad\ge
\cases{
\displaystyle 0.5 2^{l_k} {2^{a_k}\over3^{a_k}},  &\quad if
$T^{i_0+2^{l_k-a_k-1}}\omega\in
C_k$,\vspace*{2pt}\cr
\displaystyle 10^{-k'}, &\quad if $T^{i_0+2^{l_k-a_k-1}}\omega\in D_{k'}$ for some
$k<k'$,\vspace*{2pt}\cr
\displaystyle 2^{l_{k}},  &\quad if $T^{i_0+2^{l_k-a_k-1}}\omega\in C_{k'}$ for
some $k<k'$,\vspace*{2pt}\cr
\displaystyle 0.5 2^{l_k} \biggl({2\over3}\biggr)^{a_k+1},  &\quad if otherwise,}
\end{eqnarray*}
where we assumed that $l_{k'}-2 a_{k'}>l_k$ if $k'>k$.
Now
\begin{eqnarray*}
&&
|E( X_{i_0+2^{l_k-a_k-1}}|X^{i_0+2^{l_k-a_k-1}-1}_{-\infty})-
E( X_{i_0+2^{l_k-a_k-1}}|X^{i_0+2^{l_k-a_k-1}-1}_{-\infty})|\\
&&\qquad\ge{1\over4} 2^{l_k} \biggl({2\over3}\biggr)^{a_k+1}.
\end{eqnarray*}
Therefore,
\begin{eqnarray*}
{1\over N} \sum_{n=1}^{N} | E( X_{n}| X^{n-1}_0)-
E( X_{n}|
X^{n-1}_{-\infty})|&\ge&
{1\over N} {1\over4} 2^{l_k} \biggl({2\over
3}\biggr)^{a_k+1}\\
&=&
2^{-l_k+a_k} {1\over4} 2^{l_k} \biggl({2\over
3}\biggr)^{a_k+1}\\
&=&
{1\over6} \biggl({4\over3}\biggr)^{a_k}.
\end{eqnarray*}
Since $\limsup_{k\to\infty} a_k=\infty$,
the proof of Theorem \ref{counterex1} will be complete as soon as
we verify that the process is ergodic and finitarily Markovian. The first
property follows from the fact that $T$ is an ergodic transformation.
To see the second, what we need to do is to show that
the values of $f(T^{-n}\omega)$ will reveal to us more and more of the
values of
$\omega_m$ as n increases.
Almost every point is in infintely many $T^{2^{l_j-a_j}} E_j$'s. For
any such
$j$, there is a unique $i< 2^{l_j-a_j}$ such that $T^{i-2^{l_j-a_j}}
\omega\in D_j$ and this is
revealed to us by the value of $f$ at the point in the negative orbit
of $\omega$. This
information will give us the values of $\omega_m$ for all $m$ up to
$l_j-a_j$ and this completes the proof.
\end{pf}
\begin{remark}
The referee pointed out that a simpler and equivalent formulation of
the first statement of the theorem above is as follows.\vspace*{1pt}

Let $\mathcal{X}=\{0, 10^{-k}, {2^{k}\over3^{m}},
k=1,2,\ldots, m=1,2,\ldots\}$.
There exists a stationary and ergodic finitarily Markovian time series
$\{X_n\}$
taking values from $\mathcal{X}$ such that $E |X_0|<\infty$
and
\[
\limsup_{N\to\infty} {1\over N} \sum_{n=1}^{N} |E
(X_n|X^{n-1}_0)|=\infty
\]
almost surely.

[This is because $E(X_n|X^{n-1}_{-\infty})(\omega
)=E(X_0|X^{-1}_{-\infty})(T^n \omega)$ and by the
ergodic theorem
\[
\lim_{N\to\infty} {1\over N} \sum_{n=1}^{N}
E(X_n|X^{n-1}_{-\infty})=E(X_0)<\infty
\]
almost surely.]
\end{remark}
\begin{theorem} \label{counterex2}
Let $\mathcal{X}=\{0, 10^{-k}, {2^{k}\over3^{m}}
k=1,2,\ldots, m=1,2,\ldots\}$.
Suppose $h_m\dvtx \mathcal{X}^m \rightarrow\R$ is a scheme that for any
bounded ergodic
finitarily Markovian process $\{Y_n\}$
taking values from $\mathcal{X}$,
almost
surely satisfies
\[
\lim_{N\to\infty} {1\over N} \sum_{n=1}^{N}
|E(Y_n|Y^{n-1}_0)-h_n(Y_0^{n-1})|=0.
\]
Then there is an ergodic finitarily Markovian process $\{X_n\}$
taking values from $\mathcal{X}$
for which
\[
E |X_0|<\infty
\]
and
\[
\limsup_{N\to\infty} {1\over N} \sum_{n=1}^{N}
|E(X_n|X^{n-1}_0)-h_n(X_0^{n-1})|=\infty
\]
almost surely.
\end{theorem}
\begin{pf}
I. \textit{A Master process}. We shall prepare a master process
with many possibilities for constructing a
process such as in the earlier example with $l_k$ in a fashion that
will be
dictated by the estimation scheme.
For $1\le j\le n$,
define
\[
q(n,j)=(n^2+j)!
\]
and sets
\[
C_{q(n,j)}=\bigl\{\omega\dvtx \omega_i=1  \mbox{ for }1\le i<q(n,j), \omega
_{q(n,j)}=0\bigr\},
\]
clearly $P(C_{q(n,j)})=2^{-q(n,j)}$. Let
\begin{eqnarray*}
D_{q(n,j)}&=&\bigl\{\omega\dvtx \omega_i=1  \mbox{ for }1\le
i<q(n,j)-j,\\
&&\hspace*{5.1pt}\omega_{q(n,j)-j}=0,
\omega_i=1  \mbox{ for }q(n,j)-j< i<q(n,j)\bigr\}
\end{eqnarray*}
and
\[
E_{q(n,j)}=\bigcup_{i=0}^{2^{q(n,j)-j-1}-1}T^{-i} D_{q(n,j)}.
\]
Notice that
\[
E_{q(n,j)}=\bigl\{\omega\dvtx \omega_{q(n,j)-j}=0, \omega_{i}=1\mbox{,
for all }q(n,j)-j< i< q(n,j)\bigr\}
\]
and it follows that the sets $\{E_{q(n,j)}, 1 \le j \le n, n \in
\mathbf{N} \}$ are
mutually independent.
Letting
\[
u(\omega)=\sum_{n=1}^{\infty} \sum_{j=1}^{n} 10^{-q(n,j)}
I_{D_{q(n,j)}}(\omega)
\]
the master process is defined by $Y_n(\omega)=u(T^n\omega)$. For
later use, observe that the $D_{q(n,j)}$'s are disjoint.

We will need the following easy consequence of our assumption on the estimators
$h_n$, namely that for any bounded process $Y_n$ defined on $\Omega$
as in the theorem and for any
$k$ there is an integer $N_k$ and a set $H_k \subset\Omega$
with $P(H_k)\ge1-2^{-k}$ and for
all $\omega\in H_k$ and $m\ge N_k$ we have:
$|h_m(Y_0,\ldots,Y^{(k-1)})|\le{m\over10}$.\vspace*{2pt}

II. \textit{The construction}. We shall now define a
sequence $l_k$, $k=2^{a_k}+b_k$, $1\le b\le2^a$
inductively, together with functions $v_{k}$ which are bounded. As $k$
tends to infinity, the $v_k$ will converge to $v$ and we will use $u+v$
to get our desired process. We may take $v_2 =0$ to
start the inductive construction.

Assume that we have already defined $l_3<l_4<\cdots<l_{k-1}$ a
subsequence of the
$q(n,j)$'s and
\[
v_{k-1}=\sum_{i=3}^{k-1} \biggl({2^{l_i}\over3^{a_i}}\biggr) I_{C_{l_i}}
\]
we want to define $l_k$ and $v_k$.
Recalling the notation
$k-1=2^{a_{k-1}}+b_{k-1}$, we have that $b_{k-1}=b_k-1$ unless
$k-1=2^a$, in which case $a_{k-1}=a-1$ and $b_{k-1}=2^{a-1}$.

Since $v_{k-1}$ is bounded, the process defined by
\[
X_n^{(k-1)}= f_{k-1}(T^n\omega)=u(T^n\omega)+v_{k-1}(T^n\omega)
\]
is bounded.
Now, by assumption,
there is an $N_k$ and a set $H_k$ with $P(H_k)\ge1-2^{-k}$ and for
all $\omega\in H_k$ and $m\ge N_k$ we know that
\[
\bigl|h_m\bigl(X_0^{(k-1)},\ldots,X_{m-1}^{(k-1)}\bigr)\bigr|\le{m\over
10}.
\]
Choose $n$ large enough so that $2^{q(n,a_k)-a_k}>10 N_k$
and we make sure that $q(n,a_k)-a_k>10 l_{k-1}$.
Set
\[
l_k=q(n,a_k)
\]
and
\[
v_k=v_{k-1}+\biggl({2^{l_k}\over3^{a_k}}\biggr) I_{C_{l_k}}.
\]
This defines a new process
\[
X_n^{(k)}(\omega)= f_{k}(T^n\omega)=u(T^n\omega)+v_{k}(T^n\omega).
\]
It is important to observe that if for some $i_0\le2^{l_k-a_k-1}$ we have
$T^{i_0}\omega\in D_{l_k}$ then for all $0\le j\le i_0+2^{l_k-a_k-1}-1$
\[
X_j^{(k)}(\omega)=X_j^{(k-1)}(\omega).
\]
This is because the way $C_{l_k}$ is defined, we know that
$T^{i_0+2^{l_k-a_k-1}}\omega$ can be in $C_{l_k}$ which implies that earlier
iterates of $\omega$ cannot be there. Indeed,
\[
C_{q(n,j)} \subset T^{2^{l_k-a_k-1}} D_{l_k}  \qquad\mbox{for all
}q(n,j)\ge l_k,
\]
which implies that during all the later
stages of the construction the values of
$X_i^{(k-1)}$ in this range will not change. So we will have for
\[
v=\sum_{k=3}^{\infty} \biggl({2^{l_k}\over3^{a_k}}\biggr) I_{C_{l_k}}
\]
and
\[
X_n(\omega)= f(T^n\omega)=u(T^n\omega)+v(T^n\omega)
\]
that
\[
X_j(\omega)=X_j^{(k-1)} (\omega)  \qquad\mbox{for all }0\le j\le
i_0+2^{l_k-a_k-1}-1,
\]
if $T^{i_0}\omega\in D_{l_k}$.

It is clear that if the $l_k$'s are chosen large enough so that for all $k'>k$
$l_k<l_{k'}-2 a_{k'}$:
\begin{itemize}
\item the sets $\{C_k,D_k\}_{k=3}^{\infty}$ are disjoint,
\item the intervals $[l_k-a_k,l_k-1]$ are also disjoint and therefore
the sets
$E_{l_k}$ are independent.
\end{itemize}
The signaling function $u$ is bounded
and the main contributor to $f$ will be
\[
v(\omega)=\sum_{k=3}^{\infty} {2^{l_k} \over3^{a_k} }
I_{C_{l_k}}(\omega).
\]
Clearly,
\[
E(v(\omega))=\sum_{k=3}^{\infty} {2^{l_k} \over3^{a_k} } P(C_{l_k})=
\sum_{k=3}^{\infty} {1 \over3^{a_k} } =
\sum_{a=1}^{\infty} \sum_{b=1}^{2^a} {1 \over3^{a} }=
\sum_{a=1}^{\infty} \biggl({2 \over3}\biggr)^{a}<\infty.
\]
Define a process by
$
f(\omega)=u(\omega)+v(\omega)
$
and
\[
X_n(\omega)=f(T^n\omega).
\]
Note that
$X_n\in\mathcal{X}$ as advertised.

III. \textit{Checking the properties}.
Observe that $P(E_{l_k})=2^{-a_k}$ and
\[
\sum_{k=3}^{\infty} P(E_{l_k})=\sum_{a=1}^{\infty} \sum
_{b=1}^{2^a} 2^{-a}=
\sum_{a=1}^{\infty} 1= \infty.
\]
By the Borel--Cantelli lemma, a point $\omega$ belongs to $E_{l_k}$
infinitely often.
In addition, since $P(H_k)>1-2^{-k}$, almost every point will belong to
$H_k$ for
all sufficiently large $k$. Suppose then that $\omega\in
E_{l_k}\cap H_k$.
When $\omega\in E_{l_k}$,
\[
T^{i_0}\omega\in D_k  \qquad\mbox{for some }0\le i_0\le2^{l_k-a_k-1}-1.
\]
For $\omega\in E_{l_k}$, we know that $X_{i_0}(\omega)=10^{-l_k}$.
At time $i_0+2^{l_k-a_k-1}-1$,
\[
(T^{i_0+2^{l_k-a_k-1}-1}(\omega))_j =
\cases{
0,  &\quad if $j=1$,\cr
1,  &\quad if $1<j\le l_k-1$,\cr
\omega_j, &\quad  otherwise.}
\]
Let's compute for a fixed $i_0$ such that $T^{i_0}\omega\in D_{l_k}$
(i.e.,
$X_{i_0}=10^{-l_k}$)
\[
E( X_{i_0+2^{l_k-a_k}}| X^{i_0+2^{l_k-a_k}-1}_0).
\]
For $\omega\in T^{-i_0}D_{l_k}$ (i.e., $X_{i_0}=10^{-l_k}$ ) we know that
\[
(T^{i_0+2^{l_k-a_k-1}}\omega)_j =
\cases{
1,  &\quad if $1\le j \le l_k-1$,\cr
\omega_j, &\quad otherwise.}
\]
Therefore if $X_{i_0+2^{l_k-a_k-1}}>0$, then we must have
\[
T^{i_0+2^{l_k-a_k-1}} \omega\in C_{l_k}\cup\bigcup_{m>k}C_m
\cup\bigcup_{1\le n, 1\le j\le2^n\dvtx q(n,j)>l_k} D_{q(n,j)},
\]
because if $k'<k$ then $l_{k'}<l_k$ and the $C_{k'}$, are defined by zero
values of $\omega_i$ with $i<l_k$,
and similarly for $D_{q(n,j)}$ with $q(n,j) < l_k$,
\begin{eqnarray*}
&&
E( X_{i_0+2^{l_k-a_k-1}}| X^{i_0+2^{l_k-a_k-1}-1}_0)\\
&&\qquad=
{ {2^{l_k} /3^{a_k}} 2^{-l_k} +\sum_{j>k} {2^{l_j} /3^{a_j}}
2^{-l_j}+
\sum_{q(n,j)>l_k} 10^{-q(n,j)} 2^{-q(n,j)} +0 \over
P(D_{l_k}) } \\
&&\qquad\ge {({2/3})^{a_k+1} \over2\cdot2^{-l_k} }\\
&&\qquad=
{1\over2} 2^{l_k} \biggl({2\over3}\biggr)^{a_k+1}.
\end{eqnarray*}
Similarly,
\begin{eqnarray*}
&&
E( X_{i_0+2^{l_k-a_k}}| X^{i_0+2^{l_k-a_k}-1}_0)\\
&&\qquad=
{ {2^{l_k} /3^{a_k}} 2^{-l_k} +\sum_{j>k} {2^{l_j} /3^{a_j}}
2^{-l_j}+
\sum_{q(n,j)>l_k} 10^{-q(n,j)} 2^{-q(n,j)} +0 \over
P(D_{l_k}) } \\
&&\qquad\le {10^{-l_k}+\sum_{i=0}^{\infty} ({2/3}
)^{a_k+i} \over2\cdot
2^{-l_k} }\\
&&\qquad=
{1\over2} 2^{l_k} \biggl( 10^{-k-1}+\biggl({2\over3}\biggr)^{a_k}
3\biggr)\\
&&\qquad\le 4 \cdot2^{l_k} \biggl({2\over3}\biggr)^{a_k}.
\end{eqnarray*}
On the other hand, because $\omega\in H_k$ and our remark about
$X_j=X_j^{(k-1)}$ for $0\le j\le2^{l_k-a_k-1}-1$, we have that
\[
|h_{i_0+2^{l_k-a_k-1}} (X_0^{i_0+2^{l_k-a_k-1}-1})|\le
{i_0+2^{l_k-a_k-1}-1 \over
10}.
\]
Therefore, if we take $N=2^{l_k-a_k}$
\begin{eqnarray*}
{1\over N} \sum_{n=1}^{N} | E( X_{n}| X^{n-1}_0)-
h_{n} (X_0^{n-1})
|&\ge&
{1\over N} {1\over4} 2^{l_k} \biggl({2\over
3}\biggr)^{a_k+1}\\
&=&
2^{-l_k+a_k} {1\over4} 2^{l_k} \biggl({2\over
3}\biggr)^{a_k+1}\\
&=&
{1\over6} \biggl({4\over3}\biggr)^{a_k}.
\end{eqnarray*}
Since $\limsup_{k\to\infty} a_k=\infty$,
the proof of Theorem \ref{counterex2} is complete.
\end{pf}

\vspace*{-12pt}

\begin{appendix}\label{app}
\section*{Appendix}

The next result is a generalization of a result due to Elton; cf.
Theorems 2 and 4 in Elton \cite{Elton1981}.
\begin{proposition} \label{PropositionElton}
For $n=0,1,2,\ldots,$ let $\F_n$ be an increasing sequence of $\sigma
$-fields, and
$X_n$ random variables
measurable with
respect to $\F_n$, be identically distributed with
$E(|X_0|\log^+(|X_0|))<\infty$.
Let
$g_n(X_n)$ be quantizing functions so that for all $n$,
$
|g_n(X_n)-X_n|\le1,
$
and for an increasing sequence of sub $\sigma$-fields, $\G_n\subseteq
\F_n$ such
that $g_n(X_n)=Y_n$ is measurable with respect to $\G_n$, form the
sequence of
martingale differences
\[
Z_n=g_n(X_n)-E(g_n(X_n)|\G_{n-1})=Y_n-E(Y_n|\G_{n-1}).
\]
Then
%
\begin{equation}
\label{eqTHLd}
E\Biggl( \sup_{1\le n} \Biggl| {1\over n}
\sum_{i=1}^n Z_i\Biggr|\Biggr)<\infty
\end{equation}
and
%
\begin{equation}
\label{eqTHLc}
\lim_{n\rightarrow\infty}
{1\over n} \sum_{i=1}^n Z_i =0  \qquad\mbox{almost
surely.}
\end{equation}
\end{proposition}
\begin{pf}
We follow Elton \cite{Elton1981}, who gave the proof when the martingale
differences $Z_n$ are identically distributed. Write
\[
Y_n=Y'_n+Y''_n,
\]
where $|Y'_n|\le n$ and $|Y''_n|> n$. Now
\[
Z_n=Y'_n-E(Y'_n|\G_{n-1})+Y''_n-E(Y''_n|\G_{n-1}).
\]
Since for any sequence of real numbers $\{a_i\}$,
\[
\sup_{1\le n} {1\over n} \Biggl| \sum_{i=1}^n a_i\Biggr|
\le
2 \Biggl( \sup_{1\le n} \Biggl| \sum_{i=1}^n {1\over i} a_i
\Biggr|\Biggr),
\]
(cf. Lemma 7 in Elton \cite{Elton1981}), letting
\[
d_n=Y'_n-E(Y'_n|\G_{n-1})
\]
and
\[
e_n=Y''_n-E(Y''_n|\G_{n-1})
\]
we get
%
\begin{eqnarray}
&&
E\Biggl(\sup_{1\le n} {1\over n} \Biggl| \sum_{i=1}^n Z_i
\Biggr|\Biggr)\\
\label{nulladik}
&&\qquad\le  2 E\Biggl( \sup_{1\le n} \Biggl| \sum
_{i=1}^n {1\over i}
Z_i\Biggr|\Biggr)\\
\label{elso}
&&\qquad\le
2 E\Biggl( \sup_{1\le n} \Biggl| \sum_{i=1}^n {1\over i} d_i
\Biggr|\Biggr)\\
\label{masodik}
&&\qquad\quad{} +
2 E\Biggl( \sup_{1\le n} \Biggl| \sum_{i=1}^n {1\over i} e_i
\Biggr|\Biggr).
\end{eqnarray}
For (\ref{elso}) by Davis' inequality (valid for all martingale
differences cf.
e.g., Shiryayev \cite{Shiryayev1984}, page 470), we get
\[
2 E\Biggl( \sup_{1\le n} \Biggl| \sum_{i=1}^n {1\over i} d_i
\Biggr|\Biggr)\le
2 B E \Biggl[ \Biggl( \sum_{i=1}^{\infty} {1\over i^2} (d_i)^2
\Biggr)^{0.5}\Biggr]\le
2 B \Biggl[E\Biggl( \sum_{i=1}^{\infty} {1\over i^2} (d_i)^2
\Biggr)\Biggr]^{0.5}.
\]
Now,
$ E( (d_i)^2 )\le E( (Y'_i)^2 )$. But since
$|Y_i-X_i|\le1$, we get
\[
E( (Y'_i)^2 )=E\bigl( (Y_i)^2 I_{\{|Y_i|\le i\}}
\bigr)\le
E\bigl( (X_i+1)^2 I_{\{|X_i-1|\le i\}} \bigr)
\]
and the $X_i$'s are identically distributed therefore
\begin{eqnarray*}
&&
\sum_{i=1}^{\infty} {1\over i^2} E\bigl( (X_i+1)^2 I_{\{
|X_i-1|\le i\}}
\bigr)\\
&&\qquad=
\sum_{i=1}^{\infty} \Biggl( E\bigl( (X_i+1)^2 I_{\{i-1<|X_i-1|\le i\}
}\bigr)
\Biggl(\sum_{j=i}^{\infty} {1\over j^2} \Biggr)\Biggr)\\
&&\qquad\le K E( |X_0|),
\end{eqnarray*}
where $K$ is a suitable constant (cf. the last line of the proof of
Lemma 1 in
Elton \cite{Elton1981}).

For (\ref{masodik}),
\[
E|e_n|\le2 E|Y''_n|\le2 E\bigl( (1+|X_n|) I_{\{|X_n|>n-1\}}\bigr)
\]
and now $X_n$' are identically distributed. Now since $E(|X|\log
^+(|X|))<\infty$,
Lemma 2 in Elton \cite{Elton1981} implies that
\[
\sum_{n=1}^{\infty} {1\over n} E\bigl( (1+|X_n|) I_{\{|X_n|>n-1\}
}\bigr)<\infty
\]
and so
\begin{eqnarray*}
2 E\Biggl( \sup_{1\le n} \Biggl| \sum_{i=1}^n {1\over i} e_i
\Biggr|\Biggr)&\le&
2 \sum_{i=1}^n {1\over i} E|e_i|\\
&<&\infty
\end{eqnarray*}
and this completes the proof of (\ref{eqTHLd}).

Now, we prove (\ref{eqTHLc}). By (\ref{nulladik}),
\[
U_n=\sum_{i=1}^n {1\over i} Z_i
\]
is a martingale with
\[
\sup_{1\le n}
E(|U_n|)<\infty
\]
and by Doob's convergence theorem $U_n$ converges almost surely.
Then by Kronecker's lemma (cf. Shiryayev \cite{Shiryayev1984}, page 365),
\[
\lim_{n\to\infty} {1\over n} \sum_{i=1}^n Z_i=0
\]
almost surely.
The proof of Proposition \ref{PropositionElton} is complete.
\end{pf}
\begin{proposition} \label{PropositionL}
Let $\{\phi_n, \mathcal{F}_n\}$ be a martingale difference
sequence. If, for some $1<p<\infty$,
$\sup_{1\le n}E(|\phi_n|^p)<\infty$ then
%
\begin{equation}
\label{eqTHLa}
\lim_{n\rightarrow\infty} {1\over n} \sum_{i=1}^n \phi_i =0
\qquad\mbox{almost
surely}
\end{equation}
and
%
\begin{equation}
\label{eqTHLb}
E\Biggl( \sup_{1\le n} \Biggl| {1\over n} \sum_{i=1}^n \phi_i
\Biggr|^p\Biggr)<\infty.
\end{equation}
\end{proposition}
\begin{pf}
Choose a positive integer $K$ such that $K(p-1)>1$. Define
\[
f_n={1\over n} \sum_{i=1}^n \phi_i.
\]
Assume first that $1<p\le2$.
Now by Theorem 2 in von Bahr and Esseen \cite{vonBahrEsseen1965},
%
\begin{equation}
\label{eqc}
E( | f_n |^p )\le2 {n\over n^p}
\sup_{1\le i}E(|\phi_i|^p)
=2 {\sup_{1\le i}E(|\phi_i|^p)\over n^{(p-1)} }.
\end{equation}
Define
\[
F={\sum_{n=1}^{\infty}} | f_{n^K}|^p.
\]
By (\ref{eqc}), and since by assumption $\sup_{1\le n}E(|\phi
_n|^p)<\infty$,
$K(p-1)>1$,
%
\begin{equation}
\label{eqEF}
E(F)=2 \sum_{n=1}^{\infty} {\sup_{1\le i}E(|\phi_i|^p)\over
n^{K(p-1)} }<\infty.
\end{equation}
Define
\[
g_n=\max_{1\le k<(n+1)^K-n^K} | f_{n^K}-f_{n^K+k}|^p
\]
and let
\[
G=\sum_{n=1}^{\infty} g_{n}.
\]
To complete the proof of (\ref{eqTHLa}) and (\ref{eqTHLb}), it is
enough to show that
$E(F+G)<\infty$. By (\ref{eqEF}), it is enough to show that
$E(G)<\infty$. Now
for some $m=n^K+k$, $1\le k<(n+1)^K-n^K$,
\[
f_m=(f_{n^K+k}-f_{n^K})+f_{n^K}
\]
and
\[
|f_m|^p\le2^p (| f_{n^K+k}-f_{n^K}|^p +|
f_{n^K}|^p)
\le2^p ( g_n+ |f_{n^K}|^p) \le
2^p ( G+F ).
\]
Now
\[
| f_{n^K+k}-f_{n^K}|=\biggl( {1\over n^K}-{1\over
n^K+k}\biggr)
\sum_{i=1}^{n^K} \phi_i -{1\over n^K+k}\sum_{j=1}^k \phi_{n^K+j}
\]
and so
\begin{eqnarray*}
| f_{n^K+k}-f_{n^K}|^p
&\le&2^p \Biggl( \Biggl| {k\over n^K (n^K+k)}
\sum_{i=1}^{n^K} \phi_i\Biggr|^p +
\Biggl| {1\over n^K+k} \sum_{j=1}^k \phi_{n^K+j}\Biggr|^p\Biggr) \\
&\le& 2^p \Biggl( \Biggl| {(n+1)^K-n^K-1 \over n^K n^K}
\sum_{i=1}^{n^K} \phi_i\Biggr|^p +
\Biggl| {1\over n^K} \sum_{j=1}^k \phi_{n^K+j}\Biggr|^p\Biggr).
\end{eqnarray*}
Now
\[
g_n\le2^p \Biggl( \Biggl| {(n+1)^K-n^K \over n^K n^K}
\sum_{i=1}^{n^K} \phi_i\Biggr|^p +
\Biggl| {1\over n^{Kp} } \max_{1\le k<(n+1)^K-n^K} \sum_{j=1}^k
\phi_{n^K+j}\Biggr|^p\Biggr).
\]
Now by von Bahr and Eseen \cite{vonBahrEsseen1965} and Doob's
inequality (cf. e.g.,
Theorem 1, Chapter 3 in Shiryayev \cite{Shiryayev1984}),
\begin{eqnarray*}
E(g_n) &\le& 2^p \biggl( {(n+1)^K-n^K \over n^{2K} } \biggr)^p 2 n^K
\Bigl( \sup_{1\le i}E(|\phi_i|^p)\Bigr)\\
&&{}+
\biggl({p \over(p-1)}\biggr)^p {1\over n^K}
E \Biggl( \Biggl| \sum_{j=1}^{(n+1)^K-n^K} \phi_{n^K+j} \Biggr|^p
\Biggr) \\
&\le&
2^p \biggl( {(n+1)^K-n^K \over n^{2K} } \biggr)^p 2 n^K
\Bigl( \sup_{1\le i}E(|\phi_i|^p)\Bigr)\\
&&{}+ \biggl({p \over(p-1)}\biggr)^p
\biggl({ (n+1)^K-n^K \over n^K}\biggr)^p \sup_{1\le i}E(|\phi_i|^p)
\end{eqnarray*}
and the right-hand side is summable.
We have completed the proof for $1<p\le2$.

Now assume $2<p<\infty$.
By the theorem of Dharmadhikari, Fabian and Jogdeo~\cite{DFJ1968},
\[
E( | f_n |^p )\le C(p)
{\sup_{1\le i}E(|\phi_i|^p) \over n^{p/2}}.
\]
Applying this one gets that
\[
E\Biggl( \sum_{n=1}^{\infty} | f_{n}|^p\Biggr)\le
\sum_{n=1}^{\infty} C(p)
{\sup_{1\le i}E(|\phi_i|^p) \over n^{p/2}}<\infty.
\]
Thus,
\[
\sum_{n=1}^{\infty} | f_{n}|^p <\infty \qquad\mbox{almost surely}
\]
and this yields (\ref{eqTHLa}) and (\ref{eqTHLb}).
The proof of Proposition \ref{PropositionL} is complete.
\end{pf}
\begin{remark}
The referee pointed out that the second statement of the preposition
above could be proved in a simpler way as follows.
By maximal Doob inequality and Burkholder inequality, we obtain
\[
\Biggl({E\sup_n }\Biggl|{1\over n} \sum_{i=1}^n \phi_i\Biggr|^p\Biggr)^{1/ p}
\le2p \max\biggl\{1,{1\over(p-1)^2 }\biggr\}
\Biggl[ E\Biggl( \sum_{i=1}^{\infty} \biggl( {\phi_i\over i}
\biggr)^2\Biggr)^{p/2}\Biggr]^{1/ p}.
\]
Now if $p\ge2$, then by the triangle inequality in $L_{p/2}$,
\begin{eqnarray*}
\Biggl\{\Biggl[ E\Biggl( \sum_{i=1}^{\infty} \biggl( {\phi_i\over
i}\biggr)^2\Biggr)^{p/2}\Biggr]^{2/p}\Biggr\}^{1/2}
&\le&
\Biggl[ \sum_{i=1}^{\infty} \biggl(E {|\phi_i|^p\over i^p}
\biggr)^{2/ p} \Biggr]^{1/2} \\
&\le&
\Biggl( \sum_{i=1}^{\infty} {1\over i^p} \Biggr)^{1/2}\sup_i
(E |\phi_i|^p)^{1/ p}.
\end{eqnarray*}
If $p\le2$, then since
\[
\biggl(\sum_i a_i\biggr)^{p/2}\le\sum_i (a_i)^{p/2}
\]
for all positive numbers $a_i$ we get
\[
\Biggl[ E\Biggl( \sum_{i=1}^{\infty} \biggl( {\phi_i\over i}
\biggr)^2\Biggr)^{p/2}\Biggr]^{1/p}\le
\sum_{i=1}^{\infty} \biggl(E {|\phi_i|^p\over i^p} \biggr)^{1/
p} \le
\Biggl( \sum_{i=1}^{\infty} {1\over i^2} \Biggr)^{1/2}\sup_i
(E |\phi_i|^p)^{1/p}.
\]
Thus, in each case it is
\[
\Biggl(E\sup_n \Biggl| {1\over n} \sum_{i=1}^n
\phi_i\Biggr|^p\Biggr)^{1/p} \le
C_p \sup_i ( E |\phi_i|^p)^{1/p}.
\]
\end{remark}
\end{appendix}


%
\printaddresses

\end{document}